\documentstyle[12pt]{article}
\newcommand{\be}{\begin{equation}}
\newcommand{\ee}{\end{equation}}
\newcommand{\bea}{\begin{eqnarray}}
\newcommand{\eea}{\end{eqnarray}}
\newcommand{\beann}{\begin{eqnarray*}}
\newcommand{\eeann}{\end{eqnarray*}}

\newcommand{\norm}[1]{\mbox{$\|#1\|$}}

\begin{document}
\title{COMPACT PERTURBATIONS OF FREDHOLM n-TUPLES\thanks{
AMS(MOS)Subj.Classif.: 47A53, 47A55, 47A60, 47B05. Words and 
Phrases: Fredholm n-tuple, compact perturbation, index zero.}}
\author{R{\u{a}}zvan Gelca\thanks{Research partially supported by NSF.}}
\maketitle

\begin{abstract}
Let{\em T} be an operator on a Hilbert space.We show 
that the pair {\em (T,T)} can be perturbed to an invertible
pair if and only if {\em T} is Fredholm of index zero.We 
also exhibit a large class of Fredholm n-tuples acting on
a Banach space which cannot be perturbed by finite rank
operators to invertible ones.
\end{abstract}

\begin{center}
{\bf INTRODUCTION }
\end{center}
\bigskip

It is well known that a Fredholm operator of index zero on a
Banach space can be perturbed by a finite rank operator to
an invertible one.In [2] it is asked if the same property
remains true for commuting pairs of operators, or at least
if one can perturb a pair of index zero with {\em compact}
operators to get an invertible one.

There are several properties of the index that are preserved
when we pass from one operator to commuting n-tuples of
operators, for example the index is invariant under small-norm
perturbations, or under perturbations with operators in the
norm-closure of finite rank ones; in the Hilbert space case
these are exactly the compact operators [1], [4].It has 
been proved in [3] that the Koszul complex of a Fredholm 
n-tuple of index zero has a finite dimensional perturbation
to an exact complex, but the new complex is usually not a
Koszul complex of a commuting n-tuple.

In the first part of the present paper we shall prove that
on an infinite dimensional Banach space there exists a large
class of Fredholm n-tuples of index zero that cannot be
perturbed with finite rank operators to invertible ones.
Using the same idea we shall prove in the second part that
on an infinite dimensional Hilbert space no pair of the form
$(T,0)$ with {\em T} Fredholm and $ind(T)\not=0$ can be 
perturbed with compact operators to an invertible pair,
this will give a negative answer to the question raised
by Ra{\'{u}}l Curto in [2].
\bigskip

\begin{center}
{\bf 1. FINITE RANK PERTURBATIONS}
\end{center}
\bigskip

In this section we consider the case of (bounded) linear operators
acting on an infinite dimensional Banach space $ \cal X $. 
 
To each commuting n-tuple $T=(T_1,T_2,\cdots,T_n)$ of operators on
$\cal X$, we attach a complex of Banach spaces,
called the Koszul complex [5], as follows. Let 
$ \Lambda ^p= \Lambda ^p[e_1,e_2, \cdots ,e_n]$ be the p-forms on
${\bf C}^n$.Define the operator $ D_T: {\cal X} \bigotimes \Lambda^p
\rightarrow {\cal X} \bigotimes \Lambda^{p+1} $ by 
$D_T:=T_1\bigotimes E_1+T_2\bigotimes E_2+ 
\cdots +T_n\bigotimes E_n$, where $E_i\omega :=e_i\omega ,i=1,\cdots ,
n$.

The Koszul complex is

\begin{equation}
 0 \rightarrow  {\cal X} \bigotimes  \Lambda^0 
\stackrel{D_T}{\rightarrow}  
 {\cal X} \bigotimes \Lambda^1  \stackrel{D_T}{\rightarrow} 
 \cdots \stackrel{D_T}{\rightarrow} 
{\cal X} \bigotimes   
\Lambda^n \rightarrow  0.
\end{equation}

Let $H^p(T)$ be its cohomology spaces. The n-tuple $T$ is called 
invertible if $H^p(T)=0, 0\leq p\leq n$, and Fredholm if
$dim H^p(T)< \infty, 0 \leq p \leq n$, in which case we define its
index to be
$indT:=\sum_{p=0}^{n}(-1)^{p}dimH^{p}(T).$

The Taylor spectrum of $T$, denoted by $\sigma(T)$, is the set of all
$z=(z_1,z_2, \cdots ,z_n)$ in ${\bf C}^n$ such that 
$z-T=(z_1-T_1,z_2-T_2, \cdots ,z_n-T_n)$ is not invertible. It is known
that $\sigma(T)$ is a compact nonvoid set. For any holomorphic map
$f:U \rightarrow $${\bf C}^m$ on a neighborhood $U$ of $\sigma(T)$ one can
define $f(T)$ [6]; this functional calculus extends the 
polynomial calculus. By the spectral mapping theorem [6],
$f(\sigma(T))=\sigma(f(T))$.

If $T=(T_1,T_2, \cdots ,T_n)$ and 
$T'=(T_1,T_2, \cdots ,T_n,S)$ are  commuting
tuples, then we have a long exact sequence
in cohomology

$0\rightarrow  H^0(T') \rightarrow  H^0(T) 
\stackrel{\hat{S}}{\rightarrow} 
H^0(T)\rightarrow  H^1(T')
\rightarrow  H^1(T)\rightarrow  \cdots$  
\begin{equation}
 H^{p-1}(T)\rightarrow  H^p(T')
\rightarrow  H^p(T) \stackrel{\hat{S}}{\rightarrow}
 H^p(T)\rightarrow  \cdots 
\end{equation}
where ${\hat S}$ is the operator induced by $S\bigotimes 1 :  
{\cal X}
\bigotimes \Lambda^p \rightarrow {\cal X} \bigotimes \Lambda^p,
0 \leq p \leq n$. If $T'$ is invertible 
then ${\hat S}$ is an isomorphism.
As a consequence of the long exact sequence, if $T$ is Fredholm
then $T'$ is Fredholm of index zero. We shall prove that there
exists a class of (n+1)-tuples of index zero of the form
$(T_1,T_2, \cdots ,T_n,S)$ with 
$(T_1,T_2, \cdots ,T_n)$ Fredholm and
ind$(T_1,T_2, \cdots ,T_n)\not=0$ that cannot be made invertible by
finite rank perturbations.
\medskip

{\bf Lemma 1.1.} Let $(S_1,S_2,\cdots,S_n)$ 
be an invertible commuting n-tuple
and let $f:{\bf C}^n\rightarrow {\bf C}^m$ be a holomorphic function with
$f^{-1}(0)=\{0\}$. Then $f(S_1,S_2,\cdots,S_n)$ is invertible.

{\bf Proof.} Straightforward from the 
spectral mapping theorem by noticing
the fact that 
$f^{-1} (0) \bigcap
\sigma(S_1,S_2,\cdots,S_n) = \emptyset.$
\medskip

{\bf Lemma 1.2.} Let $(S_1,S_2,\cdots,S_n)$ 
be a Fredholm commuting n-tuple, with the property that
$ind (S_1,S_2,\cdots,S_n)\not=0$. 
Then there exists a sequence of
positive integers $\{m_k\}_k$ and $0\leq p_0\leq n$ 
such that
$dim H^{p_0}(S_{1}^{m_k},S_2,\cdots,S_n)
\rightarrow \infty$ for
$k\rightarrow \infty$.

{\bf Proof.} By Corollary 3.8 in [4] 
$ind(S_{1}^{m},S_2,\cdots,S_n)=
m\cdot ind(S_1,S_2,\cdots,S_n)$,  
so $dimH^p(S_{1}^{m},S_2,\cdots,
S_n)$ cannot all remain bounded.
\medskip

{\bf Theorem 1.1.} Let $(T_1,T_2,\cdots,T_n)$ be 
a Fredholm commuting n-tuple with 
$ind(T_1,T_2,$ $\cdots,T_n)\not=0$, 
and $p\in{\bf C}[z_1,z_2,
\cdots,z_n]$ with $p(0)=0$. 
Define the operator $T_{n+1}=$     
$p(T_1,T_2,\cdots,T_n)$.
Then there do not exist finite rank 
operators $R_1, R_2, \cdots,$
$R_n, R_{n+1}$ such that 
$(T_1+R_1, T_2+R_2, \cdots,$$T_{n+1}+R_{n+1})$
is an invertible comuting (n+1)-tuple.

{\bf Proof.} Suppose that such finite
 rank operators exist and let
$S_i=T_i+R_i, 1{\leq}i{\leq}n+1$.
 Applying Lemma 1.1 to the function
$f:{\bf C}^{n+1}\rightarrow 
{\bf C}^{n+1}, f(z_1, z_2, \cdots, z_n, z_{n+1})
=(z_1, z_2, \cdots, z_n, z_{n+1}-
p(z_1, z_2, \cdots, z_n))$ we get that 
$(S_1, S_2, \cdots, S_n, R)$ must 
be invertible, where 
$R=S_{n+1}-p(S_1, S_2, \cdots, S_n).$ 
Clearly, $R$ is a finite rank operator.
By applying Lemma 1.1
 to the function
$\phi : {\bf C}^{n+1} \rightarrow  
{\bf C}^{n+1},
 \phi(z_1, z_2, \cdots, z_n, z_{n+1})=
(z_1^m, z_2, \cdots, z_n, z_{n+1})$
we get that $(S_1^m, S_2, 
\cdots, S_n , R)$ is also invertible, for every
positive integer {\em m}.

Let $\{m_k\}_k$ and $p_0$ be the numbers
 obtained by applying Lemma 1.2
to the n-tuple $(S_1, S_2, \cdots, S_n)$; 
let {${\hat R}={\hat R}(m_k,p_0)$} be the operator
induced by $R$ on $H^{p_0}(S_1^{m_k}, 
S_2, \cdots, S_n)$. Because 
$(S_1^{m_k}, S_2, \cdots, S_n)$ is invertible, 
${\hat R}$ must be an isomorphism
for every $m_k$. But this is impossible since
 $dimH^{p_0}(S_1^{m_k}, S_2, \cdots, S_n)
\rightarrow  \infty$ and $rank$(${\hat R}$)
 $ \leq(^n_{p_0}) \cdot rank(R)$. 
This proves the theorem.

\bigskip
\begin{center}
{\bf 2. THE MAIN EXAMPLE}
\end{center}
\bigskip

In what follows we shall restrict 
ourselves to bounded linear
operators on an infinite dimensional 
Hilbert space $\cal H$. We shall start
with a result about the structure 
of a Fredholm operator of positive index.
\medskip

{\bf Lemma 2.1.} Let $T$ be a Fredholm 
operator with $indT>0$. 

Define 
${\cal H}_n =ker T^n{\ominus}kerT^{n-1}$. 
Then ${\cal H}_n \not=(0), n{\geq}2$.
Let $T_n=T|kerT^n$.

\begin{equation}
T_n : {\cal H}_n \bigoplus kerT^{n-1} \rightarrow  
{\cal H}_{n-1} \bigoplus kerT^{n-2} ;  T_n =
\left[
\begin{array}{clcr}
A_n & 0 \\
B_n & C_n
\end{array}
\right] .
\end{equation}
Then there exists $n_0$ such that, for $n>n_0$, $A_n$ is an
isomorphism.

{\bf Proof.} Suppose that for some $n$ ${\cal H}_n=(0)$. Then
$kerT^n=kerT^{n-1}$. Hence $kerT^{n+k}=kerT^n, \forall k \geq 0$.
But this contradicts the fact that 
$ \lim_{n \rightarrow  \infty} indT^{n+k} = \infty $.

Since ${\cal H}_n \bot kerT_{n-1}$, $T|{\cal H}_n$ is injective
and $T{\cal H}_n \bigcap kerT_{n-2} = (0)$. 
This shows that $A_n$ is injective.
But then the sequence $\{ dim {\cal H}_n \}_n$ 
is decreasing so it becomes stationary.
Let $n_0$ be such that for 
$n>n_0$, $dim{\cal H}_n = dim{\cal H}_{n-1}$. Then for $n>n_0$,
$A_n$ is an injective operator between 
finite dimensional spaces of same dimension
so it is an isomorphism.
\medskip

{\bf Lemma 2.2.} Let $T$ and ${\cal H}_n, n{\geq}2$, be as 
in the statement of previous lemma. If $S$ is an operator
that commutes with $T$, then for all $n{\geq}1$, $kerT^n$ is
an invariant subspace for $S$. Let $S_n=S|kerT^n$,
\begin{equation}
S_n : {\cal H}_n \bigoplus kerT^{n-1} \rightarrow
{\cal H}_n \bigoplus kerT^{n-1} ; S_n =
\left[
\begin{array}{clcr}
X_n & 0 \\
Y_n & Z_n
\end{array}
\right] .
\end{equation}
Then there is $n_0$ such that for $n \geq n_0$,
$X_n$ is similar to $X_{n_0}$.

{\bf Proof.} The fact that $kerT^n$ is invariant for $S$
follows from the commutativity. Let $A_n$ and $n_0$ be
as in Lemma 2.1. Then $ST=TS$ implies $S_{n-1}T_n=
T_nS_n, n \geq 2$. Therefore, $X_{n-1}A_n=A_nX_n, n
\geq 2$. For $n>n_0$ $A_n$ is an isomorphism hence
$X_n$ is similar to $X_{n-1}$. This proves the lemma.
\medskip

{\bf Lemma 2.3.} Let $(T,S)$ be an invertible commuting
pair. Then for any $n$, $S|kerT^n$ is an isomorphism of
$kerT^n$.

{\bf Proof.}  Applying Lemma 1.1 to $(T,S)$ and $f:{\bf
C}^2 \rightarrow {\bf C}^2$, $f(z_1,z_2)=(z_1^n,z_2)$
we get that $(T^n,S)$ is invertible for any $n$. By
the remarks made at the beginning of the first section, ${\hat{S}}
: H^0(T^n)  \rightarrow  H^0(T^n)$ is an isomorphism.
But $H^0(T^n)=kerT^n$, and the lemma is proved.
\medskip

{\bf Theorem 2.1.} Let $T$ be a Fredholm operator with
$indT \neq 0$. Then there do not exist compact operators
$K_1$ and $K_2$ such that $(T+K_1, K_2)$ is an invertible
commuting pair.

{\bf Proof.} Suppose such $K_1$ and $K_2$ exist. Without
loss of generality we may assume $indT>0$, otherwise we
take $T^*$ instead of $T$. We can also assume that $K_1=0$,
otherwise we can denote $T+K_1$ by $T$, and let $K_2=K$.

Consider the spaces ${\cal H}_n, n \geq 2$, obtained by
applying Lemma 2.1 to $T$, and let $K_n=K|kerT^n$.
By Lemma 2.2,
\begin{equation}
K_n : {\cal H}_n \bigoplus kerT^{n-1} \rightarrow
{\cal H}_n \bigoplus kerT^{n-1} ; K_n =
\left[
\begin{array}{clcr}
X_n & 0 \\
Y_n & Z_n
\end{array}
\right].
\end{equation}      
have the property that $X_n$ similar to $X_{n_0}$
for some $n_0$ and $n \geq n_0$. Applying Lemma 2.3
we get that the operators $X_n, n \geq 2$ are isomorphisms. 
If we denote by $r$ the spectral radius of $X_{n_0}$, then
$r>0$. From the fact that $X_n$ is similar to $X_{n_0}$
for $n \geq n_0$, (so all $X_n'$s have the same spectral radius),
it follows that ${\norm {X_n}} \geq r$.

But ${\norm{K|{\cal H}_n}}={\norm{K_n|{\cal H}_n}} \geq 
{\norm{X_n}} \geq r$ for $n \geq n_0$. Because ${\cal H}_n
\perp {\cal H}_m,n \neq m$, and ${\cal H}_n \neq (0)$
for any $n$, it follows that $K$ is not compact, a contradiction.
Therefore such $K_1$ and $K_2$ cannot exist.
\medskip

{\bf Corollary 2.1.} The pair $(T,T)$ can be perturbed by 
compacts to an invertible commuting pair if and only if
$T$ can be perturbed by a compact to an invertible operator.
\bigskip

{\bf Bibliography}

[1] Ambrozie, C.,{\em On Fredholm index in Banach spaces}, preprint, 1990;

[2] Curto, R., {\em Problems in multivariable operator
theory}, Contemporary Math., 120 (1991), 15-17;

[3] Curto, R., in {\em Surveys of some recent results 
in operator theory}, Conway, J. and Morrel, B., eds., vol. II,
Pitman Res. Notes in Math., Ser. 192, Longman Publ. Co.,
London, 1988, 25-90;

[4] Putinar, M., {\em Some invariants for semi-Fredholm systems
of essentially commuting operators}, JOT 8(1982), 65-90;

[5] Taylor, J.,L., {\em A joint spectrum for several commuting
operators}, J.Funct.Anal. 6 (1970), 172-191;

[6] Taylor, J.,L., {\em The analytic functional calculus for
several commuting operators}, Acta Math. 125 (1970), 1-38.
\medskip

Department of Mathematics, The University of Iowa, Iowa City, IA
52242.

{\em E-mail: rgelca@math.uiowa.edu}
\medskip

and

\medskip

Institute of Mathematics of Romanian Academy, P.O.Box 1-764,
70700 Bucharest, Romania.

\end{document}